\documentclass[a4,leqno,amsfonts,12pt]{article}
\usepackage[pdftex]{graphicx}
\usepackage{amssymb}
\usepackage{caption}
\usepackage{amsthm}
\newcommand{\C}{{\bf C}}

\newcommand{\R}{{\bf R}}

\newcommand{\N}{{\bf N}}

\newcommand{\de}{\delta}

\newcommand{\mb}{\mbox}

\newcommand{\beq}{\begin{equation}}
\newcommand{\eeq}{\end{equation}}

\newcommand{\ve}{\varepsilon}

\newcommand{\ov}{\overline}
\newcommand{\al}{\alpha}
\newcommand{\be}{\beta}

\newcommand{\om}{\omega}
\newcommand{\z}{\zeta}

\newcommand{\ga}{\gamma}
\newcommand{\la}{\lambda}
\newcommand{\Ga}{\Gamma}

\newtheorem{thh}{Theorem}

\newtheorem{lem}{Lemma}
\newtheorem*{deff}{Definition}
\newcommand{\La}{\Lambda}
\newcommand{\ueberschrift}{\bigskip\goodbreak\noindent\bigskip}
\newcounter{theabsatz}
\newcommand{\absatz}[1]{\stepcounter{theabsatz} \ueberschrift
	{\large \bf \arabic{theabsatz}. {#1}} \setcounter{equation}{0}}

\begin{document}
	
\begin{center}
	{\large \bf    A simple upper bound for Lebesgue constants associated with Leja
		points on the real line}
	\\[3ex] 
	
	{\bf Vladimir Andrievskii, Fedor Nazarov}\\[3ex]
	
	{\it  Department of Mathematical Sciences, Kent State University
		}\\[2ex]
	
\end{center}

\vspace*{2cm}

Running head: Lebesgue constants for Leja points on the line

\vspace*{2cm}

Mailing address:

\noindent V. Andrievskii, F. Nazarov\\   Department of Mathematical Sciences\\
Kent
State University \\ Kent, OH 44242, USA\\[4ex]

{\it E-mail addresses}:\\ andriyev@math.kent.edu, \,\, nazarov@math.kent.edu\\

{\it Phones}: (330) 672 9029,\,\, (330) 672 9013

\newpage

\begin{center}
	{\bf Abstract}
\end{center}
\vspace*{0.5cm}
	Let $K\subset \R$ be a regular compact set  and let $g(z)=g_{\ov{\C}\setminus K}(z,\infty)$ be the Green function for $\ov{\C}\setminus K$ with pole at infinity. For $\de>0$, define
$$
G(\de):=\max\{ g(z): z\in \C, \,\mb{dist}(z,K)\le 2\de\}.
$$	
Let $\{ x_n\}_{n=0}^\infty$ be a Leja sequence of points of $K$. Then the uniform norm $\|T_n\|=\La_n, n=1,2,\ldots$ of the associated interpolation operator $T_n$, i.e., the $n$-th Lebesgue constant, is bounded from above by
$$
\min_{\de>0}2n\left[\frac{\mb{diam}( K)}{\de}e^{nG(\de)}\right]^{9/8}.
$$
In particular, when $K$ is a uniformly perfect subset of $\R$, the Lebesgue constants grow at most polynomially in $n$.

To the best of our knowledge, the result is new even when $K$ is a finite union of intervals.

\vspace*{2cm}

{\it Key Words:} Leja points, Green's function, interpolation, uniformly perfect sets.

{\it AMS classification:}  41A05, 41A10

\newpage

	\absatz{Introduction}
	
Let $K\subset\R$ be any compact set. We say that a sequence of points $\{ x_n\}_{n=0}^\infty$ in $K$ is a {\it Leja sequence} if for every $k\ge 1$, $x_k$ is a point
 of global maximum of the product $\prod_{j:j<k}|x-x_j|$ on $K$ ($x_0\in K$ can be arbitrary). Let $C(K)$ be the space of continuous on $K$ functions $f:K\to \C$ endowed with the uniform norm $\|\cdot\|$ and let
 $T_n:C(K)\to C(K)$ be the {\it Lagrange interpolation operator} associated with the points $x_0,\ldots , x_{n-1}$, i.e.,
 $$
 T_nf(x)=\sum_{k=0}^{n-1} f(x_k)L_{k,n}(x)
 $$
 where 
 $$
L_{k,n}(x)=\prod_{j: 0\le j\le n-1,j\neq k}\frac{x-x_j}{x_k-x_j}
	\, .
$$
A natural question is for what classes of functions we have $T_nf\to f$
in $C(K)$ and what is the rate of convergence.  In \cite{rei, taytot, calman, pri} the reader can find  
known results,
a discussion of the rich history of this question, and further  references.

Let $E_m(f):=\inf\{\|f-p\|: p \mb{ is a polynomial of degree} \le m \}$ be the error of the best approximation of $f$ by polynomials of degree at most $m$ on $K$ and let $\La_n= \|T_n\|$ be the $n$-th {\it Lebesgue constant}.
Then, denoting by $p_{n-1}^*$ the polynomial of best approximation of degree at most $n-1$ and taking into account that $T_{n}p_{n-1}^*=p_{n-1}^*$, we can write
\begin{eqnarray*}
	\|f-T_nf\|
	&=& \|(f-p_{n-1}^*)-T_n(f-p_{n-1}^*)\|\\
	&\le& (1+\La_n)\|f-p_{n-1}^*\|=
	(1+\La_n)E_{n-1}(f).
\end{eqnarray*}
Thus, the condition $\La_n E_{n-1}(f)\to 0$ as $n\to\infty$ is sufficient for the convergence. Many natural classes of functions can be described in terms of the  rate at which $E_n(f)$ tends to $0$. For instance, let $I\supset K$ be a closed interval. Then, for every fixed $k\in \N$, the condition that the function $f$ is a restriction to $K$ of a function  having
 continuous $k$-th derivative on $I$  implies that 
$E_n(f)=O(n^{-k})$. So knowing a good upper bound for $\La_n$ and comparing it to the rate of decay of $E_n(f)$ can help one to decide quickly if the Leja interpolation scheme is guaranteed to converge.

Another reason to care about $\La_n$ is that often even if we do know a priori that $f$ is very good, the actual data we interpolate may be noisy, the noise coming either from the measurement errors, or simply from rounding in the numerical computations. In this case, we technically want to estimate not $\|f-T_nf\|$, but rather $\|f-T_n(f+\Delta f)\|$ where $\Delta f$ is the noise. Then, even when the former is small, the latter may be huge because $\|T_n\Delta f\|$ is large.
Since the noise can be completely arbitrary, it is hard to control the norm 
$\|T_n\Delta f\|$ by anything better than $\La_n \|\Delta f\|$, so, again, having a reasonably clear idea of what $\La_n$ is allows one to estimate a priori what level of noise is acceptable.

It is worth noting that in practice $\La_n$
can be evaluated a posteriori once the points $x_0,\ldots ,x_{n-1}$ have been chosen, so the a priori theoretical bounds are not necessarily a must, especially when they fall far short from what is really observed in the computations. Still, we believe that they may hold some value, if not as a prediction, then, at least, as an explanation of the pretty decent efficiency of the Leja interpolation scheme.

	\absatz{Main results}

Let $K\subset\R$ be a regular compact set, i.e., the Green function $g(z)=g_{\ov{\C}\setminus K}(z,\infty)$  of $\ov{\C}\setminus K$ with pole at infinity  is continuous on $\C$ (by definition, $g_{|K}=0$).

For $\de>0$, let
$$
G(\de):=\max\{ g(z): z\in\C,\, \mb{dist}(z,K)\le 2\de\}.
$$
\begin{thh} \label{th1}
Let $ x_0,\ldots,x_{n-1}\in K$ be a Leja sequence. Then, for every $\de>0$,
$$
\La_n\le 2n\left[ \frac{\mb{\em diam}(K)}{\de} e^{nG(\de)}\right]^{9/8}.
$$
\end{thh}
This bound, though, probably, still suboptimal for a general $K$, is fairly decent. 

Indeed, on the one hand, it shows that for any regular compact subset $K\subset\R$, $\La_n$ is subexponential. 
To see this, just take any $\ve>0$ and choose $\de>0$ so that $G(\de)<(8/9)\ve$ to obtain
$$
\limsup_{n\to \infty} \frac{\log \La_n}{n}<\ve.
$$
This already vastly extends the class of compact sets for which the subexponential upper bound was known.

On the other hand, consider a {\it uniformly perfect set} $K$. Recall that 
 according to Beardon and Pommerenke \cite{beapom} it means that
there exists a  constant $0<\ga_K<1$ such that for every $x\in K$ and every $r\in (0,\mb{diam}(K))$,
we have
$$
K\cap\{ \z\in\R:\, \ga_Kr\le|x-\z|\le r\}\neq\varnothing .
$$
Any finite union of closed intervals is uniformly perfect as well as
the classical Cantor set.
If $K\subset\R$ is uniformly perfect, then $G(\de)=O(\de^\be)$ as $\de\to 0$ for some $\be>0$ (see  \cite[pp. 562-563]{cartot}) and, choosing $\de=n^{-1/\be}$, we conclude that $\La_n=O(n^{1+(9/8)/\be})$.
In particular, when $K$ is a finite union of intervals, we have $\be=1/2$ and
$\La_n=O(n^{13/4})$.

The disadvantage of our technique is that it is unclear how to generalize it to the complex setting. The main obstacle is that the {\it  Key Lemma} (Lemma~\ref{lem2}) fails for complex points and we do not know a good substitute for it.

While Theorem \ref{th1} seems quite satisfactory from the purely theoretical point of view, in practice the issue is further complicated by the fact that one can find the maximum of a polynomial only with a certain precision, so to ensure that the Leja interpolation scheme is robust, we must also show that a small error in the maximization problem at each step does not result in high instability of the bound given by Theorem \ref{th1}. To formalize this small error possibility, we shall make the following 

 \begin{deff} Let $0<\tau\le 1$. A  sequence $x_0,\ldots,x_{n-1}\in K$ is called {\it $\tau$-quasi Leja} if for every $k=1,\ldots, n-1$, we have
$$
\prod_{j:j<k}|x_k-x_j|\ge \tau \max_{x\in K}\prod_{j:j<k}|x-x_j|.
$$
\end{deff}
In the ideal theoretical case $\tau=1$, while in practice it can be made very close to $1$ but, strictly speaking, the value $\tau=1$ is unattainable.
\begin{thh}\label{th2}
If $x_0,\ldots,x_{n-1}\in K$ is {\it $\tau$-quasi Leja}, then for every $\de>0$,
$$
 \La_n\le \frac{2}{\tau^2}n\left[ \frac{\mb{\em diam}(K)}{\tau\de} e^{nG(\de)}\right]^{9/8+2\la^{-1} \log(1/\tau)},
 $$	
where $\la=0.24565978\ldots$ is the positive root of the equation 
$e^{e^\la}(e^\la-1)=1$.
\end{thh}
Theorem \ref{th2} thus shows that the Leja scheme can exhibit at most moderate numerical instability: for $\tau$ close to $1$, the theoretical bound of Theorem \ref{th1} is just raised to some power slightly bigger than $1$.
	
	\absatz{The idea of the proof of Theorem \ref{th2}}

It is well known that
$$
\La_n=\sup_{x\in K}\sum_{k=0}^{n-1}|L_{k,n}(x)|\le n
\sup_{x\in K, 0\le k\le n-1}|L_{k,n}(x)|.
$$
Thus the main issue is to get a good bound for $|L_{k,n}(x)|$ for individual $k$. The question here is how to use the $\tau$-quasi
Leja condition in a simple but reasonably efficient way. Our suggestion is just to notice that for every $n'$ with
$k\le n'<n$, we can use the  $\tau$-quasi Leja property of $x_{n'}$ to write
\begin{eqnarray*}
\prod_{j:j<n}|x-x_j|&=&	
\prod_{j:j<n'}|x-x_j|\,\prod_{j:n'\le j<n}|x-x_j|\\
&\le&
 \tau^{-1} \prod_{j:j<n'}|x_{n'}-x_j|\,
 \prod_{j:n'\le j<n}|x-x_j|.	
\end{eqnarray*}
Repeating this trick several times, we see that if $x_0,\ldots,x_{n-1}\in K$ are
$\tau$-quasi Leja and $x=x_n\in K$ is arbitrary, then for every sequence $k=n_0<n_1<\ldots <n_m=n$, we have
$$
\prod_{j:j<n}|x-x_j|\le\tau^{-m}\prod_{j:j<k}|x_k-x_j|\,\,
\prod_{l=0}^{m-1} \,\,\prod_{j:n_l\le j<n_{l+1}}|x_{n_{l+1}}-x_j|.
$$
Now,
\begin{eqnarray*}
 |L_{k,n}(x_n)|	
	&=&
\frac{\prod_{j:0\le j\le n-1, j\neq k}|x_n-x_j|}
 {\prod_{j:0\le j\le n-1, j\neq k}|x_k-x_j|}
 =
\frac{\prod_{j:0\le j\le n-1}|x_n-x_j|}
{\prod_{j:0\le j\le n, j\neq k}|x_k-x_j|} .
		\end{eqnarray*}
\pagebreak	
	
\begin{figure}[ht]
	\includegraphics[width=1.0\textwidth]{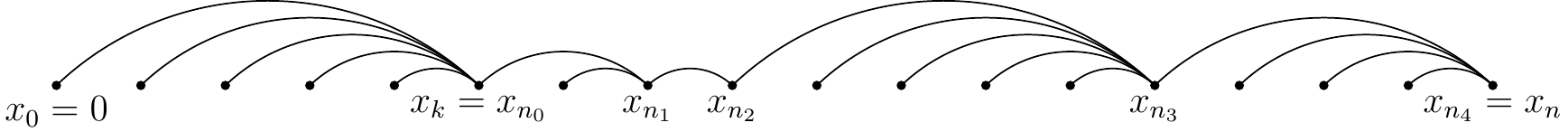}
	\caption{}
	\label{fig1}
\end{figure}

The above observations imply that we can replace the numerator on the right hand side with	
	 any product of differences similar to the one marked by arcs on Figure \ref{fig1}, in which we can change the current ``reference point" (initially $x_n$) to the last subtracted point any time at the cost of an extra $\tau^{-1}$ factor.
	 
Note that the order of points $x_j$ on this picture has nothing to do with their actual order on $\R$.

Thus, $|L_{k,n}(x_n)|\le I_\tau(x_k,\ldots,x_n)$,
where for arbitrary $x_0,\ldots,x_q\in\R$,
\beq\label{eqn1}
I_\tau(x_0,\ldots,x_q)=\min_{0=n_0<\ldots<n_m=q}\tau^{-m} \frac{\prod_{l=0}^{m-1}\prod_{j:n_l\le j<n_{l+1}}{|x_{n_{l+1}}-x_j|}}
{\prod_{j=1}^q|x_0-x_j|}\, .
\eeq
This inequality makes one tempted to forget completely about the Leja (or $\tau$-quasi Leja) property from this point on and just bound $I_\tau(x_0,\ldots,x_q)$ for an  arbitrary sequence $x_0,\ldots,x_q$ of pairwise distinct real numbers.

Unfortunately, this strategy does not quite work as it can be seen by examining the sequence $x_0=0,x_j=(-\La)^{j-1}, j\ge 1$ with $\La=1+\tau^{-1}$ being the positive root of the equation $\La^2-1=\tau^{-1}(\La+1)$.

For this sequence, it is beneficial to switch at every step, resulting in
\begin{eqnarray*}
I_\tau(x_0,\ldots,x_q)&=&	
\tau^{-q}\prod_{j=1}^q\frac{|x_j-x_{j-1}|}{|x_j-x_0|}\\	
&=& \tau^{-q}\left(1+\frac{1}{\La}\right)^{q-1}	
=\frac{1}{\tau}\left(\frac{1}{\tau}+\frac{1}{\tau+1}\right)^{q-1},
\end{eqnarray*}
which is exponential in $q$.
However, as we shall see in Section 5,
 it still does yield a universal bound
$$
I_\tau(x_0,\ldots,x_q)\le\frac{1}{\tau}\left(\frac{2}{\tau}\right)^{q-1}
$$
and the subsequent bounds
$$
|L_{k,n}(x)|\le \frac{1}{\tau}\left(\frac{2}{\tau}\right)^{n-k-1},\quad x\in K,
$$
and
$$
\La_n \le \frac{1}{\tau}\left(\frac{2}{\tau}\right)^{n}.
$$
To go below an exponential bound, we need some more information about the points $x_0,\ldots,x_{n-1}$. Fortunately, we do not need to know too much in the case $K\subset \R$ and the following simple {\it Separation Lemma} already suffices.
\begin{lem}
	\label{lem1} 
If $K$ is a regular compact set and $x_0,\ldots,x_{n-1}\in K$ is a $\tau$-quasi Leja sequence, then for every $\de>0$, 
$$
|x_i-x_j|\ge\Delta=\tau\de e^{-nG(\de)},\quad i,j=0,\ldots, n-1;i\neq j.
$$
\end{lem}
Combined with the {\it Key Lemma}  below, it immediately yields Theorem \ref{th2} and, thereby, Theorem \ref{th1} as well.

\begin{lem}
	\label{lem2} 
Let $x_0,\ldots,x_q\in\R$. Assume that $|x_0-x_j|\le D$ for all $j=1,\ldots, q$, and that  $|x_0-x_j|\ge \Delta$  for $j=1,\ldots,q-1$ with some $D\ge\Delta>0$.
Then for every $\tau\in (0,1]$, we have
$$
I_\tau(x_0,\ldots,x_q)\le\frac{2}{\tau^2}\left( \frac{D}{\Delta} \right)^{9/8+2\la^{-1}\log(1/ \tau)},
$$
where, as before, $\la$
 is the positive root of the equation 
$e^{e^\la}(e^\la-1)=1$.
\end{lem}
Note that we require no separation for $x_q$ in the assumptions of Lemma \ref{lem2}!

It remains to prove the two lemmas now.

	\absatz{Proof of the Separation Lemma}
	
Let $k$ satisfy $1\le k\le n-1$. Let $p(z)=\prod_{j:j<k}(z-x_j),\, z\in\C$.
Finally, let $M=\max_{x\in K}|p(x)|$, so $|p(x_k)|\ge\tau M.$

Consider the function $u(z)=\log(M^{-1}|p(z)|)$ on $\C\setminus K$.
It is harmonic in $\C\setminus K$ and satisfies $u(z)= 
 k\log|z|+O(1)<ng(z)$ as $z\to \infty$ and  $\limsup_{z\to K}u(z)\le 0$.
Thus, by the classical maximum principle, $u(z)\le ng(z)$ on $\C\setminus K$.
This inequality, clearly, holds on $K$ as well, so we get 
$$
|p(z)|\le Me^{ng(z)}\le Me^{nG(\de)}
$$
for every $z\in\C$ with $|z-x_k|\le 2\de$. 

By the Cauchy bound, it follows that $|p'(z)|\le M\de^{-1}e^{nG(\de)}$ when $|z-x_k|\le\de$.
In particular, when $|z-x_k|<\tau\de e^{-nG(\de)}=\Delta$, we get 
$$
|p(z)|\ge|p(x_k)|-M\de^{-1}e^{nG(\de)}|z-x_k|>\tau M-\tau M=0,
$$
so no root $x_j, j=0,\ldots,k-1$ of $p$ can lie at distance smaller than $\Delta$ from $x_k$.
This completes the proof of Lemma \ref{lem1}.

\absatz{Proof of the Key Lemma}

Replacing $x_j$ by $x_j-x_0$, we can assume  without loss of generality that
$x_0=0$. Then, replacing $x_j$ by $-x_j$ if needed, we can assume that $x_q>0$.
The ratio of the products participating in the definition of $I_\tau(x_0,\ldots,x_q)$ can be written as
$$
\frac{|X_0|}{|X_{q-1}|}\prod_{j=1}^{q-1}\frac{|X_j-x_j|}{|x_j|}\, ,
$$
where $X_j$ is the reference point at the moment of subtracting $x_j$, i.e.,
$X_j=x_{n_{l+1}}$ when $n_l\le j<n_{l+1}$. We start with $X_{q-1}=x_q$ and then, when going  from $j$ to $j-1$, can either keep the reference point (i.e., put $X_{j-1}=X_j$), or switch it to the point subtracted at the previous step (i.e., put $X_{j-1}=x_j$). The number $m$ in the prefactor $\tau^{-m}$ is just the number of switches plus $1$.

Our task is to find a good switching strategy, which, on the one hand, will allow us to control the total number of switches and, on the other hand, will keep the majority of the ratios $|X_j-x_j|/|x_j|$ small. 
The naive switching strategy is to switch every time when $|x_j|<|X_j|$.
It guarantees that for each $j=1,\ldots,q-1$ where the reference point is kept, we have
$$
\frac{|X_j-x_j|}{|x_j|}\le 2=2\frac{|X_j|}{|X_{j-1}|}
$$
and at each switch,
$$
\frac{|X_j-x_j|}{|x_j|}=\frac{|X_j-x_j|}{|X_j|}\frac{|X_j|}{|x_j|}=
\frac{|X_j-x_j|}{|X_j|}\frac{|X_j|}{|X_{j-1}|}\le 2\frac{|X_j|}{|X_{j-1}|}
$$
as well.

Thus, the full product of the ratios is at most
$$
\frac{|X_0|}{|X_{q-1}|}\prod_{j=1}^{q-1}\left(2\frac{|X_j|}{|X_{j-1}|}\right)=
2^{q-1}.
$$
Since with this strategy we may, in principle, switch at every step, the total number of switches cannot be bounded by anything better than $q-1$, so we only get the bound $I_\tau(0,x_1,\ldots,x_q)\le\tau^{-q}2^{q-1}$ mentioned earlier.

The switching strategy (or, rather, the family of switching strategies) we will consider instead is the following.  Let, as before, $\la$ be the positive root of the equation $e^{e^\la}(e^\la-1)=1$. Let $q'$ be the largest index for which $x_{q'} <0$ or $0$ if the sequence $x_1,\ldots,x_q$ contains only positive numbers.
Set $X_{q-1}=x_q$. For $j$ with $q'<j\le q-1$, put   $X_{j-1}=X_j$ if $x_j\ge e^{-\la}X_j$ and  $X_{j-1}=x_j$ if $x_j< e^{-\la}X_j$.
That is, when going over the positive tail, we switch the reference point if it becomes more than $e^\la$ times smaller after the switch and keep it otherwise.

If $q'=0$, or, which is the same, $x_j>0$ for all $j=1,\ldots, q$, then we always have
$$
\frac{|X_j-x_j|}{|x_j|}\le\frac{|X_j|}{|X_{j-1}|}\, ,\quad 1\le j\le q-1.
$$
Indeed, if we made no switch, then the right hand side is $1$, while the left hand one is either $(x_j-X_j)/x_j<1$ if $x_j\ge X_j$, or
$$
\frac{X_j-x_j}{x_j}\le\frac{(1-e^{-\la})X_j}{e^{-\la}X_j}= e^\la-1<1
$$
if $e^{-\la}X_j\le x_j<X_j$.

On the other hand, if we switch, then $|x_j|=|X_{j-1}|$ and, since $0<x_j<e^{-\la}X_j$,
we have
$$
|X_j-x_j|=X_j-x_j<X_j=|X_j|,
$$ 
so the inequality holds again. 

Also, since each switch makes the reference point $e^\la$ times smaller, we cannot have more than $\la^{-1}\log(D/\Delta)$ switches. 
This is clear if $x_q\ge\Delta$. But if $x_q<\Delta$, then we cannot have any switches at all, so
the estimate is still valid.

Thus, in this case
$$
I_\tau(0,x_1,\ldots,x_q)\le\tau^{-(\la^{-1}\log(D/\Delta)+1)}
\frac{|X_0|}{|X_{q-1}|}\prod_{j=1}^{q-1}\frac{|X_j|}{|X_{j-1}|}=
\frac{1}{\tau}\left(\frac{D}{\Delta}\right)^{\la^{-1}\log(1/\tau)}\, .
$$
Otherwise, for $j=q'$, keep both options (switching to $x_{q'}$  and staying with $X_{q'}$) available. Denote the corresponding reference points by $b_1=X_{q'}>0$ and $-a_1=x_{q'}<0$.

From this step on, we shall always have two options for the current reference point with the switching rules as follows. When going left from $q'-1$, if $x_j\not\in (-e^{-\la}a_1, e^{-\la} b_1)$, just keep the reference point as it was (either $-a_1$ or $b_1$, whichever was chosen as $X_{q'-1}$).

However, if $x_j\in (-e^{-\la}a_1, e^{-\la} b_1)$ and $x_j>0$, then force the switch from $b_1$ to $b_2=x_j$ and allow the switch from $-a_1$ to $b_2$, so, beyond this step, we shall have the options $b_2$ and $-a_1$ instead of $b_1$ and $-a_1$. For notational convenience, we will denote them $-a_2$ and $b_2$, $a_2$ being just the same as $a_1$.
Similarly, if $x_j<0$, then $x_j$ becomes $-a_2$, $b_2$ stays the same as $b_1$, and the switch from $-a_1$ to $-a_2$ is forced while the switch from $b_1$ to $-a_2$ is optional.

 Now repeat the same procedure with the interval 
$(-e^{-\la}a_2, e^{-\la} b_2)$ instead of $(-e^{-\la}a_1, e^{-\la} b_1)$, and so on until we reach $j=0$ with two options $-a_l$ or $b_l$ for $X_0$ with some $l>0$.
Our first task will be to bound the total number of switches.
Consider first the case when $x_q\ge\Delta$. Then
as long as we go along the positive tail, each switch decreases $X_j$ at least $e^\la$ times, so the number of switches made during this part of the strategy is at most
$$
\la^{-1}\log\frac{X_{q-1}}{X_{q'}}\le \la^{-1}\log\frac{D}{X_{q'}}.
$$
Then there may be one switch to $x_{q'}$. Beyond that, every time we allow a switch, the product $a_sb_s$ decreases at least $e^\la$ times, so the remaining number of switches is at most
$$
\la^{-1}\log\frac{a_1b_1}{a_lb_l}\le\la^{-1}\log\frac{DX_{q'}}{\Delta^2}\, .
$$
Adding everything up, we conclude that in this case we can have at most $2\la^{-1}
\log(D/\Delta)+1$ switches, so
\beq\label{eq1}
\tau^{-m}\le\tau^{-(2\la^{-1}\log ( D/\Delta)+2)}=\frac{1}{\tau^2}
\left(\frac{D}{\Delta}\right)^{2\la^{-1}\log(1/\tau)}\, .
\eeq
Now suppose that $x_q<\Delta$. Then we have no switches in the positive tail and $b_1=X_{q'}=x_q$.
We still may have one switch to $x_{q'}$. After that we cannot decrease $b_1$, so every time we switch, we decrease $a_s$ at least $e^\la$ times and the remaining number of switches is at most $\la^{-1}\log(a_1/a_l)\le\la^{-1}\log(D/\Delta)$.

Thus, in this case we can have at most $\la^{-1}\log(D/\Delta)+1$ switches and
$$
\tau^{-m}\le\tau^{-(\la^{-1}\log ( D/\Delta)+2)}=\frac{1}{\tau^2}
\left(\frac{D}{\Delta}\right)^{\la^{-1}\log(1/\tau)}\, ,
$$
which is an even stronger bound than that in (\ref{eq1}).

Now it is time to estimate
$$
\frac{|X_0|}{|X_{q-1}|}\prod_{j=1}^{q-1}\frac{|X_j-x_j|}{|x_j|}\, .
$$
The same argument as in the positive case shows that
$$
\prod_{j:j>q'}\frac{|X_j-x_j|}{|x_j|}\le\frac{X_{q-1}}{X_{q'}}=\frac{X_{q-1}}{b_1} .
$$
For $j=q'$, we have $\frac{|X_j-x_j|}{|x_j|}=\frac{a_1+b_1}{a_1}.$ Thus
\beq\label{eqq}
\frac{1}{|X_{q-1}|}\prod_{j:j\ge q'}\frac{|X_j-x_j|}{|x_j|}\le
\frac{a_1+b_1}{a_1b_1}.
\eeq
This leaves us with $|X_0|\prod_{j=1}^{q'-1}(|X_j-x_j|/|x_j|)$.
The product here is not unique: we have a whole family of admissible strategies, not a single one. So we shall estimate some multiplicative average of this quantity over all of them.

Let $q_2>\ldots>q_l$ be the indices at which the values of $a_s$ and $b_s$ change, i.e., $x_{q_{s+1}}\in (-e^{-\la}a_s, e^{-\la} b_s)$ while 
 $x_j\not\in (-e^{-\la}a_s, e^{-\la} b_s)$ for $q_{s+1}<j<q_s, s=1,\ldots ,l$
 (we set $q_1=q', q_{l+1}=0$ here).
 Let 
 $$
 \al_s=\frac{b_s}{a_s+b_s}\, ,\quad  \be_s=\frac{a_s}{a_s+b_s}\, .
 $$
 Then for every $j$ with $q_{s+1}<j< q_s$, we have
 $$
 \left(\frac{|-a_s-x_j|}{|x_j|}	\right)^{\al_s}
\left(\frac{|b_s-x_j|}{|x_j|}	\right)^{\be_s}\le 1
$$
(see Elementary Inequality 1 in the Appendix).

Thus, if we denote
$$
\overline P_s=\prod_{j:q_{s+1}<j<q_s}\frac{|-a_s-x_j|}{|x_j|},\quad
\overline Q_s=\prod_{j:q_{s+1}<j<q_s}\frac{|b_s-x_j|}{|x_j|},
$$ 
we have $\overline P_s^{\al_s}\overline Q_s^{\be_s}\le 1$ as well.

Consider now $j=q_{s+1}$. If $s<l$, then two cases are possible: $x_j=-a_{s+1}$ or $x_j=b_{s+1}$.

Assume that $x_j=-a_{s+1}$. Then, since $X_j$ is either $-a_s$ or $b_s$, we have $|X_j-x_j|/|x_j|$ equal to
either $(a_s-a_{s+1})/a_{s+1}$ or  $(b_s+a_{s+1})/a_{s+1}$.
By  Elementary Inequality 2,
\begin{eqnarray*}
\left[\frac{a_s-a_{s+1}}{a_{s+1}}\right]^{\al_s}
\left[\frac{b_s+a_{s+1}}{a_{s+1}}\right]^{\be_s}&\le&
\frac{a_s}{a_{s+1}}\left[\frac{\min(a_s,b_s)}{\min(a_{s+1},b_s)}\right]^{1/8}\\
&=&\frac{a_sb_s}{a_{s+1}b_{s+1}}\left[\frac{\min(a_s,b_s)}{\min(a_{s+1},b_{s+1})}\right]^{1/8}
\end{eqnarray*}
(recall that in this case $b_{s+1}=b_s$). The case $x_j=b_{s+1}$ is symmetric to the considered one and results in the same bound.

Thus, putting
$$
P_s=\overline P_s\frac{a_s-a_{s+1}}{a_{s+1}},\quad
Q_s=\overline Q_s\frac{b_s+a_{s+1}}{a_{s+1}},
$$
we get
$$
P_s^{\al_s}Q_s^{\be_s}\le
\frac{a_sb_s}{a_{s+1}b_{s+1}}\left[\frac{\min(a_s,b_s)}{\min(a_{s+1},b_{s+1})}
\right]^{1/8} ,\quad s=1,\ldots, l-1.
$$
For $s=l$, we just need to add $|X_0|$ to the product $\prod_{j=1}^{q_l-1}(|X_j-x_j|/|x_j|)$,
which results in either $P_l=\overline P_l a_l$ or $Q_l=\overline Q_l b_l$.
Thus in this case
$P_l^{\al_l}Q_l^{\be_l}\le a_l^{\al_l}b_l^{\be_l}.$
Multiplying these
 estimates out, we get
\begin{eqnarray}
\prod_{s=1}^l(P_s^{\al_s}Q_s^{\be_s})&\le&
\frac{a_1b_1}{a_{l}b_{l}}\left[\frac{\min(a_1,b_1)}{\min(a_{l},b_l)}\right]^{1/8}
a_l^{\al_l}b_l^{\be_l}\nonumber\\
\label{eq2}
&=&
\frac{a_1b_1}{a_{l}^{\be_l}b_{l}^{\al_l}}\left[\frac{\min(a_1,b_1)}{\min(a_{l},b_l)}\right]^{1/8}.
\end{eqnarray}
Now observe that for every admissible strategy under our rules, the product
$$
|X_0|\prod_{j=1}^{q'-1}\frac{|X_j-x_j|}{|x_j|}
$$
equals to $\prod_{s=1}^l R_s$ where each $R_s$ is either $P_s$ or $Q_s$ and we can start with both $P_1$ and $Q_1$, after which we can always follow $P_s$ by $P_{s+1}$ and 
$Q_s$ by $Q_{s+1}$, but we also can follow $Q_s$ by $P_{s+1}$ if $a_{s+1}<a_s, b_{s+1}=b_s$ and $P_s$ by $Q_{s+1}$ if $a_{s+1}=a_s, b_{s+1}<b_s$.
\newpage
In other words, possible products correspond to the paths on a diagram like the one on Figure 2.

\begin{figure}[ht]
	\includegraphics[width=1.0\linewidth]{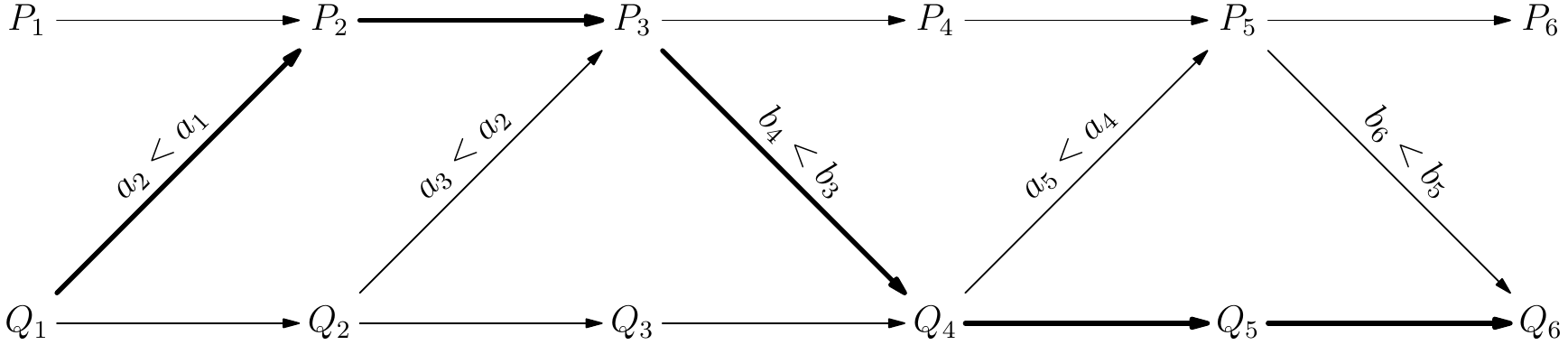}
	\caption{}
	\label{fig2}
\end{figure}

For a path $\pi$ on the diagram, let $\Ga_\pi$ be the corresponding product of $P'$s and $Q'$s.
For instance, if $\pi$ is the path following the thick edges on Figure 2, then  $\Ga_\pi =Q_1P_2P_3Q_4Q_5Q_6$.

We want to show that a certain multiplicative average $\prod_\pi\Ga_\pi^{\om_\pi}$ with some \newline $\om_\pi>0,\sum_\pi \om_\pi=1$ is exactly $\prod_{s=1}^l (P_s^{\al_s}Q_s^{\be_s})$ 
in the sense that if $P_s$ and $Q_s$ on the diagram are viewed as free variables, then the equality
$$
\prod_\pi\Ga_\pi^{\om_\pi}=\prod_{s=1}^l(P_s^{\al_s} Q_s^{\be_s})
$$
becomes an algebraic identity.
 We prove it by induction on the length $l$ of the diagram.

If $l=1$, then the only paths are $P_1$ and $Q_1$, so we can just put $\om_{P_1}=\al_1,
\om_{Q_1}=\be_1$.

Suppose now that the statement holds for $l-1$. Note that if we remove $P_1$ and $Q_1$ from a diagram of length $l$ (together with three edges coming out of them), we shall get a diagram of the same kind but of length $l-1$ and starting with $P_2$ and $Q_2$.
Thus, by the induction assumption, we can find weights $\om_{\pi'},\om_{\pi''}>0$ with
$\sum_{\pi'}\om_{\pi'}+
\sum_{\pi''}\om_{\pi''}=1$,
where $\pi'$ and $\pi''$ run over all paths starting with $P_2$ and $Q_2$ respectively, such that
$$
\prod_{\pi'}\Ga_{\pi'}^{\om_{\pi'}}\prod_{\pi''}\Ga_{\pi''}^{\om_{\pi''}}
=\prod_{s=2}^l (P_s^{\al_s}Q_s^{\be_s}).
$$
Comparing the powers at $P_2$ and $Q_2$, we see that we must have $\sum_{\pi'}\om_{\pi'}=\al_2$ and $\sum_{\pi''}\om_{\pi''}=\be_2$.
Now if $a_2<a_1$, say (i.e., if the switch from $Q_1$ to $P_2$ is allowed), we have $\al_2>\al_1$ and the admissible paths $\pi$ in the full diagram are $P_1\pi', Q_1\pi'$ and $Q_1\pi''$. Put
$$
\om_{P_1\pi'}=\frac{\al_1}{\al_2}\om_{\pi'},\quad
\om_{Q_1\pi'}=\left( 1-\frac{\al_1}{\al_2}\right)\om_{\pi'},\quad
\om_{Q_1\pi''}=\om_{\pi''}.
$$
Then in the product $\prod_\pi\Ga_\pi^{\om_\pi}$ for the full diagram, the power of $P_1$ is
$$
\frac{\al_1}{\al_2}\sum_{\pi'}\om_{\pi'}=\frac{\al_1}{\al_2}\al_2=\al_1,
$$
the power of $Q_1$ is
$$
\left( 1-\frac{\al_1}{\al_2}\right)\sum_{\pi'}\om_{\pi'}
+\sum_{\pi''}\om_{\pi''}=
\left( 1-\frac{\al_1}{\al_2}\right)\al_2+\be_2=\al_2+\be_2-\al_1=
1-\al_1=\be_1,
$$
and the powers of $P_s$ and $Q_s$ with $s\ge 2$ are exactly the same as in the product
$\prod_{\pi'}\Ga_{\pi'}^{\om_{\pi'}}\prod_{\pi''}\Ga_{\pi''}^{\om_{\pi''}}$, i.e., $\al_s$ and $\be_s$ respectively.

The case when $b_2<b_1$ (so $\be_2>\be_1$ and the switch from $P_1$ to $Q_2$ is possible) is similar just with the roles of $P'$s and $Q'$s as well as $\al'$s and $\be'$s swapped.

The upshot is that combining (\ref{eqq}) and (\ref{eq2}), we conclude that              
there exists a switching strategy in our family such that
\beq\label{eq3}
\frac{|X_0|}{|X_{q-1}|}\prod_{j=1}^{q-1}\frac{|X_j-x_j|}{|x_j|}\le
\frac{a_1+b_1}{a_{l}^{\be_l}b_{l}^{\al_l}}\left[\frac{\min(a_1,b_1)}{\min(a_{l},b_l)}
\right]^{1/8}.
\eeq
Now consider two cases.

{\bf Case 1:} $x_q\ge\Delta$. Then $a_l,b_l\ge\Delta$ and the right hand side of (\ref{eq3}) is at most
$$
\frac{2D}{\Delta}\left(\frac{D}{\Delta}\right)^{1/8}=
2\left(\frac{D}{\Delta}\right)^{9/8}.
$$
{\bf Case 2:} $x_q<\Delta$. Then $b_1=b_l=x_q<\Delta$, so the second factor on the right hand side equals $1$, while the first one can be estimated as
$$
\frac{a_1+b_1}{a_l^{\be_l}b_l^{\al_l}}=
\frac{a_l+b_l}{a_l^{\be_l}b_l^{\al_l}}
\frac{a_1+b_l}{a_l+b_l}
\le 2 \frac{a_1+b_l}{a_l+b_l}
\le 2\frac{a_1}{a_l}\le 2\frac{D}{\Delta}
$$
(see Elementary Inequality $3$).

Thus, in both cases, the right hand side of (\ref{eq3}) does not exceed $2(D/\Delta)^{9/8}$.
Combining this estimate with the previously obtained bound (\ref{eq1}), we get the statement of the lemma.

\absatz{Appendix: Elementary inequalities}

\noindent {\bf Elementary Inequality $1$}. Let $a,b>0, x\not\in(-e^{-\la}a,e^{-\la}b)$. Then
$$
\left[\frac{|x+a|}{|x|}\right]^{b/(a+b)}
\left[\frac{|x-b|}{|x|}\right]^{a/(a+b)}\le 1.
$$
Since the statement does not change if we multiply $x,a$ and $b$ by any positive number, we can assume that $a+b=1$. Our inequality then becomes
$$
\left[\frac{|x+a|}{|x|}\right]^{b}
\left[\frac{|x-b|}{|x|}\right]^{a}\le 1.
$$
 If $x\not\in [-a,b]$, we can just use the concavity of the function
 $t\mapsto\log t$ for $t>0$ to write
 $$
\left[\frac{|x+a|}{|x|}\right]^{b}
\left[\frac{|x-b|}{|x|}\right]^{a}
\le b\frac{|x+a|}{|x|}+a\frac{|x-b|}{|x|}=1.
$$
Suppose now that $x\in [-a,b]\setminus (-e^{-\la}a,e^{-\la}b)=
[-a,-e^{-\la}a]\cup[e^{-\la}b,b]$.

Since the inequality does not change if we replace the triple $(a,b,x)$ by $(b,a,-x)$, we can assume without loss of generality that $x\in[-a,-e^{-\la}a]$, i.e.,
$x=-ae^{-\Lambda}$ with $0\le\Lambda\le \la.$

Then we want to prove that
$$
\left[\frac{a(1-e^{-\Lambda})}{ae^{-\Lambda}}\right]^b
\left[\frac{b+ae^{-\Lambda}}{ae^{-\Lambda}}\right]^a
=[e^\Lambda-1]^b\left[ 1+\frac{b}{a}e^\Lambda\right]^a\le 1.
$$
However, since $1+t\le e^t$ for every $t\ge 0$,
this follows from
$$
[e^\Lambda-1]^be^{(b/a)e^\Lambda a}=
[(e^\Lambda-1)e^{e^\Lambda}]^b\le
[(e^\lambda-1)e^{e^\lambda}]^b=1.
$$

\noindent {\bf Elementary Inequality $2$}.
Let $A,B>0,0<a<A$. Then
$$
\left[\frac{A-a}{a}\right]^{B/(A+B)}
\left[\frac{B+a}{a}\right]^{A/(A+B)}\le 
\frac{A}{a}
\left[\frac{\min(A,B)}{\min(a,B)}\right]^{1/8}.
$$
Since the statement does not change  if we multiply $A,B,a$ by any positive number, we can assume without loss of generality that $A=1$. Multiplying both sides by $a$, we see that we need to prove that
$$
(1-a)^{B/(1+B)}(B+a)^{1/(1+B)}\le
\left[\frac{\min(1,B)}{\min(a,B)}\right]^{1/8}.
$$
Since 
$$
(1-a)^{B/(1+B)}(1+aB)^{1/(1+B)}\le
\frac{B}{1+B}(1-a) +\frac{1}{1+B}(1+aB)=1
$$
by the concavity of logarithm, it is enough to show that
$$
\left[\frac{B+a}{1+aB}\right]^{1/(1+B)}\le
\left[\frac{\min(1,B)}{\min(a,B)}\right]^{1/8}.
$$
If $B\le 1$, then
$
1+aB-(B+a)=(1-B)(1-a)\ge 0$,
so the left hand side is at most $1$ and the right hand side is at least $1$.

Otherwise $\min(1,B)=1,\min(a,B)=a$ and we arrive at
$$
\left[\frac{B+a}{1+aB}\right]^{1/(1+B)}\le\left[\frac{1}{a}\right]^{1/8},\quad
0<a<1, B>1.
$$
Put $a=1-t, 0<t<1$, and rewrite the inequality as
$$
\frac{1}{1+B}\left[\log\left(1-\frac{t}{B+1}\right)
-\log\left(1-\frac{Bt}{B+1}\right)
\right]\le\frac{1}{8}\log\frac{1}{1-t}\, ,
$$
or
$$
\frac{1}{1+B}\left[\varphi\left(\frac{Bt}{B+1}\right)
-\varphi\left(\frac{t}{B+1}\right)
\right]\le\frac{1}{8}\left[\varphi(t)-\varphi(0)\right],
$$
where $\varphi(t)=-\log(1-t)$. The left hand side is then $\frac{B-1}{(1+B)^2}t$
times the average of $\varphi'$ over the interval $\left[\frac{t}{B+1}, \frac{Bt}{B+1}
\right]$, while the right hand side is $t$ times the average of $\varphi'$ over the interval $[0,t]$. Since the intervals are concentric and $\varphi'(t)=1/(1-t)$ is convex, the latter average is greater, so it suffices to show that
$$
\frac{B-1}{(B+1)^2}\le\frac{1}{8}\, .
$$
However $(B+1)^2=(B-1+2)^2\ge(2\sqrt{2(B-1)})^2=8(B-1)$
and we are done.

\noindent {\bf Elementary Inequality $3$}.
Let $a,b>0$. Then
$$
\frac{a+b}{a^{\frac{a}{a+b}}
	b^{\frac{b}{a+b}}}\le 2.	
$$
Again, the inequality is invariant under the multiplication of $a,b$ by any positive number, so we can assume that $a+b=1$. Then $b=1-a$ and we are to prove that
$$
\frac{1}{a^a(1-a)^{1-a}}\le 2.
$$
However $a\mapsto a\log a +(1-a)\log (1-a)$ is a convex function on $(0,1)$ symmetric around $a=1/2$,
so its minimum is attained at $1/2$ and equals $\log(1/2)$.

\noindent {\bf Elementary Inequality $4$}.
The positive root $\la$ of the equation 
$e^{e^\la}(e^\la-1)=1$ satisfies $\la>1/5$.

 Indeed, since the function is increasing, it suffices to prove that
$e^{e^{1/5}}(e^{1/5}-1)<1$.
Note that $e^{-1/5}>1-1/5=4/5,$ so $e^{1/5}<5/4$.
Thus, it is enough to show that $e^{5/4}<4$. Since $e<3$, we may check that $3^{5/4}<4$ instead. However $3^5=243<256=4^4$ and the result follows.

	\noindent{\bf Acknowledgements.}
	The authors are  grateful to  N. Levenberg, I. Pritsker, L. Reichel, and V. Totik
	for their helpful comments.

\end{document}